\documentclass[a4paper]{article}

\usepackage{amsmath}
\usepackage{amssymb}
\usepackage{amsthm}
\usepackage{mathrsfs}

\theoremstyle{plain}
\newtheorem{thm}{Theorem}[section]
\newtheorem{prop}[thm]{Proposition}
\newtheorem{lem}[thm]{Lemma}

\theoremstyle{definition}
\newtheorem{defn}[thm]{Definition}
\newtheorem{eg}[thm]{Example}

\theoremstyle{remark}
\newtheorem{rem}[thm]{Remark}

\numberwithin{equation}{section}

\DeclareMathOperator{\dom}{dom}		
\DeclareMathOperator{\tr}{tr}		

\def \bQ {\mathbb Q}		
\def \bR {\mathbb R}		

\def \cB {\mathcal B}		%
\def \cC {\mathcal C}		%
\def \cD {\mathcal D}		%
\def \cE {\mathcal E}		%
\def \cH {\mathcal H}		%
\def \cK {\mathcal K}		%
\def \cL {\mathcal L}		%
\def \cP {\mathcal P}		%
\def \cT {\mathcal T}		%
\def \sA {\mathscr A}		%
\def \sB {\mathscr B}		%
\def \sF {\mathscr F}		%
\def \sG {\mathscr G}		%
\def \sI {\mathscr I}			%
\def \ve {\vec e}
\def \vh {\vec h}

\def \vv {\vec v}

\def \vy {\vec y}

\def \fB {\mathbf B}
\def \ff {\mathbf f}
\def \fg {\mathbf g}
\def \fU {\mathbf U}
\def \fV {\mathbf V}

\def \EH {\cE(H)}
\def \EXi {\cE(\Xi)}
\def \Hpi {\cH_\pi}
\def \Hp {\cH_1}
\def \Hm {\cH_{-1}}

\def \llangle {\langle\!\langle}
\def \rrangle {\rangle\!\rangle}

\begin{document}

\pagestyle{plain}
\bibliographystyle{plain}

\title{A central limit theorem for stochastic heat equations \\in random environment}
\author{Lu \textsc{Xu}\footnote{l-xu@math.kyushu-u.ac.jp} \quad \textit{Kyushu University}}
\date{}
\maketitle

\noindent\textbf{Abstract} In this article, we investigate the asymptotic behavior of the solution to a one-dimensional stochastic heat equation with random nonlinear term generated by a stationary, ergodic random field. We extend the well-known central limit theorem for finite-dimensional diffusions in random environment to this infinite-dimensional setting. Due to our result, a central limit theorem in $L^1$ sense with respect to the randomness of the environment holds under a diffusive time scaling. The limit distribution is a centered Gaussian law whose covariance operator is explicitly described. It concentrates only on the space of constant functions. 





\section{Introduction and main result}

Homogenization of finite-dimensional diffusions in stationary and ergodic environments is a well-studied topic, including periodic and quasiperiodic environments as special cases. In the early works \cite{K79} and \cite{PV81}, central limit theorems are established for diffusions driven by random, self-adjoint operators of divergence type. Diffusions without drift are considered in \cite{PV82}. In \cite{O83}, an invariance principle in a quenched sense is obtained for diffusions in ergodic, almost surely $C^2$-smooth environment, through a study on the fundamental solutions corresponding to its generator. 

A good review of finite-dimensional results can be found in \cite[Sect. 9]{KLO12}, where the main approach is to record the environment viewed from the particle as a Markov process. Every sample path of the diffusion is decomposed into the sum of an additive functional of the environment process and a martingale with stationary and ergodic increments. The general method to prove a central limit theorem for additive functional of Markov process is firstly developed in \cite{KV86} for reversible case and be extended later in \cite{O98}, \cite{OS95} and \cite{V95} to non-reversible case where a sector condition holds. We rely on these results in this article, but in an infinite-dimensional setting. 

Our aim is to extend this strategy to an infinite-dimensional, nonlinear system. We study the homogenization of the solution to a stochastic partial differential equation in random environment. The equation considered here is a stochastic heat equation on the unit interval $[0, 1]$ driven by standard space--time white noise and having a random nonlinear term. Different from the finite-dimensional model, the law of the nonlinear term is supposed to be stationary and ergodic under only constant shifts, that is, a group of transformations indexed by $\bR$. We adopt this setting since the Laplacian in the equation is preserved only by these transformations, which is necessary for obtaining the Markov property of the environment process. 

The nonlinear term is supposed to be the composition of a gradient-type part and a divergence-free part. If only deterministic gradient-type environment is adopted, this equation is used to describe the motion of a flexible Brownian string in some potential field, see \cite{F83}. Moreover, if the environment degenerates to a periodic nonlinear term, the model is closely related to the dynamical sine-Gordon equation (see, e.g., \cite{HS16}), and an invariance principle is obtained in \cite{X17}. The divergence-free part can be added since it preserves the equilibrium state. For homogenization of finite-dimensional diffusions in divergence-free random field, we refer to \cite{O88}. 

To state our model, we first introduce some basic notations. Throughout this article let $H = L^2[0, 1]$. The inner product and the norm on $H$ are denoted by $\langle \cdot, \cdot \rangle_H$ and $\|\cdot\|_H$. Let $E$ be the Banach space $C[0, 1]$ equipped with the uniform topology. Denote by $E_0$ the subspace of $E$ consisting of functions which vanish at $0$. Denote by $\mu_0$ the Wiener measure on $E_0$ induced by a standard Brownian motion. Since the sample path of a Brownian motion is almost surely H\"older continuous with order less than $\frac{1}{2}$, we fix some $\alpha \in (0, \frac{1}{2})$ and introduce $E^\alpha$ as the space of $\alpha$-H\"older continuous functions on $[0, 1]$. $E$, $E^\alpha$ and $E_0$ are treated as subspaces of $H$. 

Recall that a real-valued function $f$ defined on $H$ is said to be Fr\'echet differentiable, if at every $h \in H$ there exists some $Df(h) \in H$ such that 
$$\lim_{h' \in H; \|h'\|_H \downarrow 0} \frac{f(h + h') - f'(h) - \langle Df(h), h' \rangle_H}{\|h'\|_H} = 0. $$
$Df$ is called the Fr\'echet derivative of $f$. Higher-order derivatives are defined inductively. For a positive integer $k$, by $C_b^k(H)$ we denote the class of all bounded functions on $H$ which have bounded and continuous Fr\'echet derivatives up to the $k$-th order. 

Now we state our model precisely. Suppose that $(\Sigma, \sA, \bQ)$ and $(\Omega, \sF, \mathbb P)$ are two complete probability spaces, the latter of which is equipped with a filtration of $\sigma$-fields $\{\sF_t \subseteq \sF; t \geq 0\}$ satisfying the usual conditions. Let $W(t, x)$ be a standard cylindrical Brownian motion defined on $(\Omega, \sF, \mathbb P)$ and adapted to $\sF_t$, and 
$$\{(V(\sigma, u), B(\sigma, u)) \in \bR \times H; (\sigma, u) \in \Sigma \times H\} $$
be an $\bR \times H$-valued random field over $H$ on $(\Sigma, \sA, \bQ)$. Suppose that $V$ is Fr\'echet differentiable in $u$ for almost all $\sigma$, and let $U$ be the $H$-valued random field defined by $U = DV + B$. For fixed $\sigma \in \Sigma$, consider a one-dimensional stochastic heat equation with homogeneous Neumann boundary conditions and initial condition $v \in E$. 
\begin{equation}\label{spde}\left\{ 
\begin{aligned}
&\partial_t u(t, x) = \frac{1}{2}\partial_x^2 u(t, x) - U(\sigma, u(t)) + \dot W(t, x), &&t > 0, x \in (0, 1); \\
&\partial_x u(t, x)|_{x = 0} = \partial_x u(t, x)|_{x = 1} = 0, &&t > 0; \\
&u(0, x) = v(x), &&x \in [0, 1]. 
\end{aligned}\right. 
\end{equation}
Equation \eqref{spde} is called a stochastic heat equation in random environment $U$, and every $\sigma \in \Sigma$ is called a fixed environment. To make sure that for fixed $\sigma$, \eqref{spde} has a strong solution in the space of continuous functions, assume that \vspace{0.5em}\\
\textbf{(A1)} $U(\sigma, \cdot) = DV(\sigma, \cdot) + B(\sigma, \cdot)$ is a bounded and Lipschitz continuous map from $H$ to $H$ for $\bQ$-almost all $\sigma \in \Sigma$. \vspace{0.5em}\\
Under \textbf{(A1)}, for almost all $\sigma$, the unique solution $u^{\sigma, v}(t, x)$ to \eqref{spde} is continuous in $t$ and $\alpha$-H\"older continuous in $x$ (see \cite{F83, F91}) for $\alpha < \frac{1}{2}$. Therefore, $\{u^{\sigma, v}(t); t \geq 0\}$ forms a continuous Markov process taking values in $E^\alpha$. Consider the stochastic process on the product space $(\Sigma \times \Omega, \sA \otimes \sF)$ defined as 
$$u^v(t) \triangleq u^{\cdot, v}(t), \ t \geq 0. $$
It is called the solution to \eqref{spde} in random environment, which is the main object we want to study in this article. Sometimes we omit the initial condition $v$ and write it as $u(t)$ for short. To simplify the discussion on the generator of the environment process defined in Definition \ref{environment process} later, also assume that \vspace{0.5em}\\
\textbf{(A2)} $\sup_{\Sigma \times H}\left\{|V| + \|DV\|_H + \|B\|_H\right\} < \infty$. \vspace{0.5em}\\
To continue, we discuss the path space and the distribution of $\{V, B\}$. Let 
$$C^{1, 0}(H; \bR \times H) \triangleq \left\{\phi = (v, b): H \rightarrow \bR \times H~\left|~v \in C^1(H; \bR), b \in C(H; H)\right.\right\}. $$
Due to the regularity of $u^{\sigma, v}(t, \cdot)$, the distribution of $u^{\sigma, v}(t)$ depends only on the law of the subfield $\{\tilde V, \tilde B\} \triangleq \{V(u), B(u)\}_{u \in E^\alpha}$. Hence, let 
$$\Sigma_{\text{path}}^0 \triangleq \left\{\tilde \phi: E_\alpha \rightarrow \bR \times H~\left|~\exists \phi \in C^{1, 0}(H; \bR \times H), \text{ s.t. } \tilde \phi = \phi|_{E^\alpha}\right.\right\}. $$
Denote by $\Sigma_{\text{path}}$ the completion of $\Sigma_{\text{path}}^0$ under the Fr\'echet metric 
\begin{equation}\label{frechet metric}
d_\Sigma(\sigma_1, \sigma_2) = \sum_{k = 1}^\infty \frac{1}{2^k} \cdot \frac{d_{\Sigma, k}(\sigma_1, \sigma_2)}{1 + d_{\Sigma, k}(\sigma_1, \sigma_2)}, 
\end{equation}
where for $\sigma_1 = (\tilde v_1, \tilde b_1)$ and $\sigma_2 = (\tilde v_2, \tilde b_2) \in \Sigma$, 
$$d_{\Sigma, k}(\sigma_1, \sigma_2) = \sup_{u \in E_k^\alpha}\left\{|\tilde v_1(u) - \tilde v_2(u)| + \|D\tilde v_1(u) - D\tilde v_2(u)\|_H + \|\tilde b_1(u) - \tilde b_2(u)\|_H\right\}, $$
$$E_k^\alpha = \left\{u \in E^\alpha, |u(x)| \leq k, \frac{|u(x) - u(y)|}{|x-y|^\alpha} \leq k, \forall x, y \in [0, 1]\right\}. $$
Since every $E_k^\alpha$ is compact, $(\Sigma_{\text{path}}, d_\Sigma)$ is a Polish space under the metric $d_\Sigma$. Equip it with the Borel $\sigma$-field and adopt it as the path space of $\{\tilde V, \tilde B\}$. For $c \in \bR$, let $\tau_c$ be the transformation on $\Sigma_{\text{path}}$ defined by 
\begin{equation}\label{shift operators}
\tau_c \circ \tilde \phi = \tilde \phi(\cdot + c\mathbf 1), \ \forall \tilde \phi \in \Sigma_{\text{path}}, 
\end{equation}
where $\mathbf1$ stands for the constant function $\mathbf 1(x) \equiv 1$ on $[0, 1]$. Since $\tau_{c_1} \circ \tau_{c_2} = \tau_{c_1 + c_2}$, $\{\tau_c\}$ forms a group. Given a probability measure $p$ on $\Sigma_{\text{path}}$, we say $p$ is ergodic if for measurable set $A \subseteq \Sigma_{\text{path}}$ such that $p(A \Delta \tau_c[A]) = 0$ for all $c \in \bR$, we have either $p(A) = 1$ or $p(A) = 0$, where $\Delta$ stands for the symmetric difference. Let $P_{\tilde V, \tilde B}$ be the distribution of $\{\tilde V, \tilde B\}$ on $\Sigma_{\text{path}}$ and assume that \vspace{0.5em}\\
\textbf{(A3)} $P_{\tilde V, \tilde B}$ is stationary and ergodic under $\{\tau_c; c \in \bR\}$, and 
\begin{equation}\label{non-degenerate}
P_{\tilde V, \tilde B}(\tau_c\tilde\phi = \tilde\phi, \ \forall c \in \bR) < 1. 
\end{equation}
The translation invariance and ergodicity of the environment are usually necessary for studying homogenization in such a model. Since $\{\tilde\phi~|~\tau_c\tilde\phi = \tilde\phi, \ \forall c \in \bR\}$ is obviously a translation invariant set, \eqref{non-degenerate} is equivalent to say that $P_{\tilde V, \tilde B}(\tau_c\tilde\phi = \tilde\phi, \ \forall c \in \bR) = 0$. If \eqref{non-degenerate} fails to hold, $U(\sigma, u)$ is periodic in the sense $U(\sigma, u) = U(\sigma, u + c\mathbf1)$. For a stochastic heat equation with local, periodic and gradient-type nonlinear term, an invariance principle is obtained in \cite{X17}. Finally, we assume that \vspace{0.5em}\\
\textbf{(A4)} $\exists$ a measurable set $\Sigma_1$ such that $\tau_c[\Sigma_1] = \Sigma_1$, $\bQ(\Sigma_1) =1$ and for all $\sigma \in \Sigma_1$, 
\begin{equation}
\int_{E_0} e^{-2V(\sigma, u)}\langle Df(u), B(\sigma, u) \rangle_H\mu_0(du) = 0 
\end{equation}
holds for all $f$ on $H$ such that $f(u) = f^\dagger(\langle u, \varphi_1 \rangle_H, \ldots, \langle u, \varphi_M \rangle_H)$ for some $M \geq 1$, $f^\dagger \in C_b^1(\bR^M)$ and $\varphi_1, \ldots, \varphi_M \in H$. \vspace{0.5em}\\
\textbf{(A4)} is equivalent to say that $\delta(e^{-2V}B) = 0$ holds for all $\sigma \in \Sigma_1$, where $\delta$ is the divergence operator adjoint to the Malliavin derivative (see, e.g., \cite[p. 35, Definition 1.3.1]{N95}). One can see in Lemma \ref{stationarity-frozen} and Proposition \ref{stationarity} that a non-gradient field satisfying \textbf{(A4)} preserves the equilibrium state in both fixed and random environment. Our main result, a central limit theorem for $u^v(t)$, is stated as follows. 

\begin{thm} \label{clt}
Under \textbf{(A1)} to \textbf{(A4)}, $u(t)/\sqrt{t}$ satisfies the central limit theorem in $L^1$ sense with respect to the environment and the limit distribution concentrates on the space of constant functions, i.e., for any bounded continuous function $f$ on $E$, 
\begin{equation} \label{clt-equation}
\lim_{t \rightarrow \infty}E_{\bQ} \left|E_{\mathbb P} \left[f\left(\frac{u(t)}{\sqrt{t}}\right)\right] - \int_{\mathbb{R}} f(y\mathbf1)\Phi_a(y)dy\right| = 0, 
\end{equation}
where $y\mathbf1$ is the function on $[0, 1]$ taking constant value $y \in \bR$, and $\Phi_a$ is the probability density function of a 1-d centered Gaussian law with variance $a^2$ defined in \eqref{diffusion constant} below. Furthermore there exists some strictly positive constant $C$ depending only on $V$ such that $C \leq a^2 \leq 1$. 
\end{thm}

\begin{rem}
A question naturally rises up is whether central limit theorem holds $\bQ$-almost surely. We would like to point out here, that such quenched result usually relies on dimensions (cf. \cite{B11}). The most challenging part in proving a quenched result is to show the sublinearity of the corrector, which usually involves methods like heat kernel estimates or Sobolev inequality. Since these methods are not expectable in our setting, the quenched version of Theorem \ref{clt} cannot be achieved easily. 
\end{rem}

In Sect. 2, we introduce the evolution of the environment seen from the solution, which is recorded as an ergodic Markov process with values in a new space consisting of series of environments. In Sect. 3, we prove Theorem \ref{clt} mainly based on a sector condition and some general arguments in \cite[Sects. 2, 9]{KLO12}. Before these contents, we close this section with two examples, showing that both periodic and quasiperiodic nonlinear terms are included in our model as special cases. 

\begin{eg}
Take $\Sigma = [0, 1]$, $\sA = \cB(\Sigma)$ and $\bQ$ to be the Lebesgue measure. Suppose $V$ to be a measurable function on $[0, 1] \times \bR$ such that 
$$V(x, \cdot) \in C^1(\bR), \ V(x, y) = V(x, y+1), \ \forall x \in [0, 1], \ y \in \bR. $$
Define the random field $(V, B)$ for all $\sigma \in \Sigma$ and $u \in E$ as 
$$V(\sigma, u) = \int_0^1 V(x, u(x) + \sigma)dx $$
and $B(\sigma, u) \equiv 0$. Assume that both $V$ and $\frac{d}{dy}V$ are uniformly bounded, then \textbf{(A1)} to \textbf{(A4)} are fulfilled. Denote by $u^{\sigma, v}$ the solution to \eqref{spde}, then $u^{\sigma, v - \sigma}(t, x) + \sigma$ solves the stochastic heat equation with the periodic nonlinear term $\frac{d}{dy}V(x, u(t, x))$. 
\end{eg}

\begin{eg}
For $d \geq 2$ take $\Sigma = [0, 1]^d$, $\sA = \cB(\Sigma)$ and $\bQ$ to be the Lebesgue measure. Suppose $V$ to be a measurable function on $[0, 1] \times \bR^d$ such that 
$$V(x, \cdot) \in C^1(\bR^d), \ \forall x \in [0, 1]; $$
$$V(x, \vy) = V(x, \vy + \ve_i), \ \forall x \in [0, 1], \ \forall \vy \in \bR^d, \ \forall i = 1, \ldots, d; $$
$$\sup_{[0, 1] \times \bR^d} \left\{|V|, |\partial_{y_1}V|, \ldots, |\partial_{y_d}V|\right\} < \infty, $$
where $\ve_i$ is the unit vector in $\bR^d$ of $i$-th coordinate. Let $(\lambda_1, \ldots, \lambda_d) \in \bR^d$ be a vector with rationally independent coordinates. Define the random field $(V, B)$ for all $\sigma \in \Sigma$ and $u \in E$ as 
$$V(\sigma, u) = \int_0^1 V(x, \lambda_1u(x) + \sigma_1, \ldots, \lambda_du(x) + \sigma_d)dx $$
and $B(\sigma, u) \equiv 0$ for all $\sigma \in \Sigma$ and $u \in E$. It gives the quasiperiodic model. 
\end{eg}

\section{The environment process}

Since the statement in Theorem \ref{clt} depends only on the law of $(V, B)$, without loss of generality we may and will fix $\Sigma$ to be the Polish metric space $(\Sigma_{\text{path}}, d_\Sigma)$ described in Sect. 1. We take $\sA$ to be the Borel $\sigma$-algebra on $\Sigma$ associated to $d_\Sigma$, and let $\bQ = P_{\tilde V, \tilde B}$. As defined in \eqref{shift operators}, there is a group of transformations $\{\tau_c; c \in \bR\}$ on $\Sigma$ such that $\bQ$ is stationary, ergodic and satisfies that 
\begin{equation}\label{strong non-degenerate}
\bQ(\Sigma_0) = 0, 
\end{equation}
where $\Sigma_0 \triangleq \{\sigma \in \Sigma~|~\tau_c\sigma = \sigma, \forall c \in \bR\}$ is a subset of $\Sigma$. 

This section is devoted to the construction of the environment process associated to \eqref{spde}. To apply the strategy in \cite[Sect. 9]{KLO12}, the environment process needs to be a Markov process taking values in a Polish metric space, having a stationary and ergodic probability measure. Given an environment $\sigma \in \Sigma$ and some $v \in E$, to shift $\sigma$ with $v$ it is natural to consider a map $\xi: [0, 1] \rightarrow \Sigma$ such that 
$$\xi(x) = \tau_{[v(x)]}\sigma, \ \forall x \in [0, 1]. $$
In the above expression, since $\xi(\cdot) = \tau_{[v(\cdot) - v(0)]}[\xi(0)]$ where $v(\cdot) - v(0) \in E_0 = \{u \in E; u(0) = 0\}$, $\xi$ belongs to the set $\Xi$ defined as follows: 
\begin{equation}\label{environment set}
\Xi \triangleq \{\xi: [0, 1] \rightarrow \Sigma~|~\exists (\sigma, v) \in \Sigma \times E_0 \text{ s.t. } \xi(\cdot) = \tau_{[v(\cdot)]}\sigma\}. 
\end{equation}
Every $\xi \in \Xi$ is called \textit{a series of environments}, which describes the environments seen from a function. To construct a metric on $\Xi$, we need the following lemma. 

\begin{lem}\label{key lemma}
For each $\xi$ such that $\xi(0) \notin \Sigma_0$ in \eqref{strong non-degenerate}, there is only one pair of $\sigma \in \Sigma$ and $v \in E_0$ satisfying that $\xi(\cdot) = \tau_{[v(\cdot)]}\sigma$. 
\end{lem}

\begin{proof}
For $\sigma \in \Sigma$, define $\ker(\sigma) \triangleq \{c \in \bR~|~\tau_c\sigma = \sigma\}$, then $\Sigma_0 = \{\sigma \in \Sigma~|~\ker(\sigma) = \bR\}$. Define a subset $\Xi_0 \subseteq \Xi$ to be 
\begin{equation}\label{xi 0}
\Xi_0 = \{\xi \in \Xi~|~\xi(0) \in \Sigma_0\}. 
\end{equation}
For $\xi \notin \Xi_0$, suppose that $\xi = \tau_{[v(\cdot)]}\sigma = \tau_{[v'(\cdot)]}\sigma'$ for some $(\sigma, v)$ and $(\sigma', v') \in \Sigma \times E_0$. Since $v(0) = v'(0) = 0$, it is obvious that $\sigma = \sigma' = \xi(0)$. Notice that $\ker(\sigma)$ does not contain any non-degenerate interval $[a, b]$, otherwise the group property of $\{\tau_c\}$ would imply that $\ker(\sigma) = \bR$ and then $\xi \in \Xi_0$. However, simple calculation shows that 
$$\tau_{[v(x) - v'(x)]}\sigma = \sigma, \ \forall x \in [0, 1]. $$
Hence, the image of $v - v'$ contains no non-degenerate interval. As $v - v'$ is continuous on $[0, 1]$ and $v(0) - v'(0) = 0$, this yields that $v = v'$. 
\end{proof}

In view of Lemma \ref{key lemma}, one is able to define a one-to-one map $\xi \mapsto (\sigma_\xi, v_\xi)$ from $\Xi$ to $\Sigma \times E_0$ such that $\xi(\cdot) = \tau_{[v_\xi(\cdot)]}\sigma_\xi$. Indeed $\sigma_\xi = \xi(0)$, and when $\xi \notin \Xi_0$, $v_\xi$ is uniquely determined due to Lemma \ref{key lemma}. For those $\xi \in \Xi_0$, we can simply take $v_\xi(x) \equiv 0$. Equip $\Xi$ with the metric $d_\Xi$ defined as 
\begin{equation}\label{metric}
d_\Xi(\xi, \eta) = d_\Sigma(\xi(0), \eta(0)) + \mathrm{sup}_{\, x \in [0, 1]}|v_\xi(x) - v_\eta(x)|, \ \forall \xi, \eta \in \Xi, 
\end{equation}
where $d_\Sigma(\cdot, \cdot)$ is the Fr\'echet metric in \eqref{frechet metric}. Observing that $(\Xi, d_\Xi)$ may not be complete, we take the completion and still denote it by $\Xi$. Since both $(E_0, |\cdot|_{\mathrm{sup}})$ and $(\Sigma, d_\Sigma)$ are Polish spaces, also is $(\Xi, d_\Xi)$. Equip $\Xi$ with the Borel $\sigma$-algebra $\sG$, then $\xi \mapsto (\xi(0), v_\xi)$ becomes a measurable map with respect to $\sA \otimes \sB(E_0)$. Now the environment process can be defined as follows. 

\begin{defn}\label{environment process}
The environment process associated with \eqref{spde} is defined as 
$$\xi_t = \xi_t(\cdot) \in \Xi, \text{ s.t. } \xi_t(x) = \tau_{[u^{\sigma, v}(t, x)]}\sigma, \ \forall x \in [0, 1]. $$
\end{defn}

Recall that $(\Omega, \sF, \mathbb P)$ is the probability space where the noise $\dot W(t, x)$ in \eqref{spde} is defined. $\{\xi_t\}_{t \geq 0}$ defined in Definition \ref{environment process} forms a $\Xi$-valued stochastic process on the product space $(\Sigma \times \Omega, \sA \otimes \sF)$. To express the nonlinear term in \eqref{spde} by $\xi_t$, we introduce the derivative operator on $\Xi$. For $\ff: \Xi \rightarrow \bR$, let $f^\sigma$ be the function defined as $f^\sigma(v) = \ff(\tau_{[v(\cdot)]}\sigma)$ for all $v \in E$. 

\begin{defn}\label{derivative}
A function $\ff$ on $\Xi$ is called differentiable, if for every $\sigma$, $f^\sigma$ can be extended to a Fr\'echet differentiable function $\bar f^\sigma$ on $H$, and its derivative is defined as $\cD\ff(\xi) = Df^\sigma(v)$ for $\xi = \tau_{[v(\cdot)]}\sigma$. Higher-order derivatives are defined inductively. In particular, for differentiable $\ff$, $\cD\ff(\xi) = \cD\ff(\xi, \cdot) \in H$ and for twice differentiable $\ff$, $\cD^2\ff(\xi)$ is a bounded bilinear functional on $H \times H$. 
\end{defn}

The nonlinear term $DV(\sigma, u^{\sigma, v}(t)) + B(\sigma, u^{\sigma, v}(t))$ can be replaced by a function of $\xi_t$. Indeed, recall that every $\sigma = (v, b) \in \Sigma = \Sigma_{\text{path}}$ can be continuously extended to an element in $C^{1, 0}(H; \bR \times H)$, we can define 
$$\left(\fV(\xi), \fB(\xi)\right) = [\xi(0)](v_{\xi}) \in \bR \times H, \ \forall \xi \in \Xi, $$
where $v_\xi \in E_0$ is defined by the one-to-one map from $\Xi$ to $\Sigma \times E_0$ mentioned above. By Definition \ref{derivative} and the definition of $\Sigma_{\text{path}}$, it is not hard to see that 
$$\cD\fV(\xi_t) + \fB(\xi_t) \overset{d.}= DV(\sigma, u^{\sigma, v}(t)) + B(\sigma, u^{\sigma, v}(t)) = U(\sigma, u^{\sigma, v}(t)). $$

In summary, without loss of generality we can assume the following framework. Suppose $(\Sigma, d_\Sigma)$ to be some Polish metric space equipped with the Borel $\sigma$-field $\sA$ and some Borel probability measure $\bQ$. A group of measurable transformations $\{\tau_c; c \in \bR\}$ is defined on $\Sigma$ such that $\bQ$ is stationary, ergodic and satisfies \eqref{strong non-degenerate}. Define $\Xi$ as in \eqref{environment set} and equip it with the metric $d_\Xi$ in \eqref{metric} and the Borel $\sigma$-algebra $\sG$. It suffices to prove Theorem \ref{clt} for the following equation instead of \eqref{spde}: 
\begin{equation}\tag{\ref{spde}'}\label{spde'}\left\{ 
\begin{aligned}
&\partial_t u^{\sigma, v}(t, x) = \frac{1}{2}\partial_x^2 u^{\sigma, v}(t, x) - \cD\fV(\xi_t) - \fB(\xi_t) + \dot W(t, x), &&t > 0, x \in (0, 1); \\
&\partial_x u^{\sigma, v}(t, x)|_{x = 0} = \partial_x u^{\sigma, v}(t, x)|_{x = 1} = 0, &&t > 0; \\
&u^{\sigma, v}(0, x) = v(x), &&x \in [0, 1], 
\end{aligned}\right. 
\end{equation}
where $(\fV, \fB)$ is a random variable on $(\Xi, \sG)$ taking values in $\bR \times H$. The random field generated by $(V(\sigma, u), B(\sigma, u)) = (\fV(\tau_{[u(\cdot)]}\sigma), \fB(\tau_{[u(\cdot)]}\sigma))$ is supposed to satisfy \textbf{(A1)}, \textbf{(A2)} and \textbf{(A4)} in Sect. 1. 

Before stating the Markov property of $\xi_t$, we prepare some notations. Let $\cB_b(E)$ and $\cB_b(\Xi)$ be the collections of all bounded measurable functions on $E$ and $\Xi$, respectively. For fixed $\sigma$, $u^{\sigma, v}(t)$ defines a continuous Markov process with values in $E^\alpha \subseteq E$. Denote by $\cP^\sigma_t$ the associated Markov semigroup defined on $\cB_b(E)$. Denote by $\mathbb P_v^\sigma$ the law of $u^{\sigma, v}$ for fixed $\sigma$ and $v$, and by $\mathbb E_v^\sigma$ the corresponding expectation. With these notations, we have the next proposition. 

\begin{prop}\label{markov property}
The environment process $\{\xi_t, t \geq 0\}$ is a Markov process on $\Xi$. Let $\{\cP_t, t \geq 0\}$ be the Markov semigroup on $\cB_b(\Xi)$ determined by $\xi_t$, then 
\begin{equation}\label{semigroup}
\cP_t\ff(\xi) = \cP_t^{\xi(0)}f^{\xi(0)}(v_\xi), \ \forall \ff \in \cB_b(\Xi), 
\end{equation}
where $f^\sigma$ is the function on $E$ defined by $f^\sigma(v) = \ff(\tau_{[v(\cdot)]}\sigma)$. 
\end{prop}

\begin{proof}
First notice that if $\xi_0 \in \Xi_0$ in \eqref{xi 0} then the environment $\sigma \in \Sigma_0$. This furthermore implies that $\xi_t = \xi_0$ for all $t \geq 0$ and $f^{\xi(0)}$ is a constant function on $E$. Hence, the Markov property and \eqref{semigroup} hold obviously in this case. 

Pick $\sigma \notin \Sigma_0$, $v \in E$ and let $\xi = \tau_{[v(\cdot)]}\sigma$. Let $P^\sigma(t; v, \cdot)$ be the transition probability of the solution $u^{\sigma, v}(t, \cdot)$ in a fixed environment. Since $u^{\tau_c\sigma, v} = u^{\sigma, v + c} - c$, 
$$P^\sigma(t; v, \cdot) = P^{\tau_c\sigma}(t; v - c, \cdot - c), \ \forall c \in \bR. $$
Since $\xi \notin \Xi_0$, the discussion below Lemma \ref{key lemma} yields that $v_\xi = v(\cdot) - v(0)$, thus by taking $c = v(0)$ in the equation above, we obtain that 
\begin{equation}\label{shift probability}
P^{\sigma}(t; v, \cdot) = P^{\xi(0)}(t; v_\xi, \cdot - v(0)), 
\end{equation}
Pick some $t > 0$, $h > 0$, $n \geq 0$, $0 \leq t_1 \leq \ldots \leq t_n < t + h$ and bounded measurable functions $\ff_1, \ldots, \ff_n, \fg$ on $\Xi$ arbitrarily. We have 
$$\mathbb E_v^\sigma \left[\prod_{j=1}^n\ff_j(\xi_{t_j})\fg(\xi_{t + h}) \right] = \mathbb E_v^\sigma \left[\prod_{j=1}^nf_j^\sigma(u^{\sigma, v}(t_j))g^\sigma(u^{\sigma, v}(t + h))\right]. $$
Using the Markov property of $u^{\sigma, v}(t)$ for fixed $\sigma$, the expectation above equals to 
$$\mathbb E_v^\sigma \left[\prod_{j=1}^nf_j^\sigma(u^{\sigma, v}(t_j))\int_E g^\sigma(u)P^\sigma(h; u^{\sigma, v}(t), du)\right]. $$
Applying \eqref{shift probability} with $\xi = \xi_t$ and $v = u^{\sigma, v}(t)$, this expectation equals to 
\begin{equation*}
\begin{aligned}
&\mathbb E_v^\sigma \left[\prod_{j=1}^n\ff_j(\xi_{t_j})\int_E g^\sigma(u)P^{\xi_t(0)}(h; v_{\xi_t}, d(u - u^{\sigma, v}(t, 0)))\right] \\
=~&\mathbb E_v^\sigma \left[\prod_{j=1}^n\ff_j(\xi_{t_j})\int_E g^{\xi_t(0)}(u)P^{\xi_t(0)}(h; v_{\xi_t}, du)\right]. 
\end{aligned}
\end{equation*}
Therefore, $\xi_t$ is a Markov process and 
$$\cP_t\ff(\xi) \triangleq \int_E f^{\xi(0)}(u)P^{\xi(0)}(t; v_\xi, du) = \cP_t^{\xi(0)}f^{\xi(0)}(v_\xi), \ \forall t \geq 0 $$
is the associated Markov semigroup. 
\end{proof}

Next we construct the stationary and ergodic measure of $\{\cP_t, t \geq 0\}$. Before discussing the environment process, we give a lemma concerning the (infinite) stationary measure in frozen environment $\sigma \in \Sigma$. Notice that in \cite[Sect. 4]{F83}, the same problem is discussed for the case that $B(\sigma, \cdot) = 0$ and $V(\sigma, \cdot)$ is local. 

\begin{lem}\label{stationarity-frozen}
For almost all $\sigma \in \Sigma$, $\{\cP_t^\sigma; t \geq 0\}$ admits a stationary measure 
$$\pi^\sigma(dv) = \exp(-2V(\sigma, v))\mu(dv), $$
where $(E, \mu)$ is the infinite measure determined by a stochastic process $\{w(x); x \in [0, 1]\}$ satisfying that $w(\cdot) - w(0)$ is a one-dimensional Brownian motion, while $w(0)$ subjects to the Lebesgue measure on $\bR$. Moreover if $\fB \equiv 0$, then $\pi^\sigma$ is reversible. 
\end{lem}

This lemma follows from a Galerkin approach. For $k \geq 1$, let $\lambda_k = 2(k - 1)^2\pi^2$ be the eigenvalues of $-\frac{1}{2}\partial_x^2$ on $H$ with Neumann boundary condition. The corresponding eigenvectors are $h_1 = \mathbf1$ and $h_k(x) = \sqrt 2\cos[2(k - 1)\pi x]$ for $x \in [0, 1]$ and $k \geq 2$. With these notations, we sketch the proof for completeness. 

\begin{proof}[Proof of Lemma \ref{stationarity-frozen}]
In view of \textbf{(A4)}, since $(V, B)(\sigma, v + c\mathbf1) = (V, B)(\tau_c\sigma, v)$ holds for all $c \in \bR$, it is easy to see that if $\sigma$ belongs to $\Sigma_1$, 
\begin{equation}\label{A4'}
E_{\pi^\sigma} \left[\langle Df, B(\sigma, \cdot) \rangle_H\right] = \int_\bR E_{\mu_0} \left[e^{-2V(\tau_c\sigma, \cdot)}\langle Df, B(\tau_c\sigma, \cdot)\right]dc = 0, 
\end{equation}
for $f(v) = f^\dagger(\langle v, h_1 \rangle_H, \ldots, \langle v, h_N \rangle_H)$ with $f^\dagger \in C_0^1(\bR^N)$. Since here only fixed environment is considered, we temporarily omit $\sigma$ and write $V$, $B$, $u^v(t)$ in short of $V(\sigma, \cdot)$, $B(\sigma, \cdot)$ and $u^{\sigma, v}(t)$. Let $\nu(dv) = e^{-2V(v)}\mu(dv)$, and we prove that $\nu(dv)$ is stationary for $\{u^v(t); t \geq 0\}$ under \eqref{A4'}. 

First assume that $V$ and $B$ are finite-dimensional dependent. Precisely, for some $N \geq 1$, $V^\dagger \in C_b^1(\bR^N)$ and $B^\dagger = (B_1^\dagger, \ldots, B_N^\dagger) \in C_b(\bR^N; \bR^N)$, we have 
\begin{equation}\label{cylindrical V and B}
V(v) = V^\dagger(\langle v, h_1 \rangle_H, \ldots, \langle v, h_N \rangle_H),~B(v) = \sum_{k=1}^N B_k^\dagger(\langle v, h_1 \rangle_H, \ldots, \langle v, h_N \rangle_H)h_k. 
\end{equation}
Let $X_t^k = \langle u^v(t), h_k \rangle_H$, then $X_t = (X_t^1, \ldots, X_t^N) \in \bR^N$ solves 
$$dX_t^k = -\lambda_kX_t^kdt - \partial_kV^\dagger(X_t) - B_k^\dagger(X_t)dt + d\langle W(t, \cdot), h_k \rangle_H, \ X_0^k = \langle v, h_k \rangle_H. $$
Define $\Psi_k(x) = \sqrt{\lambda_k\pi^{-1}}e^{-\lambda_kx^2}$ for $k \geq 1$ and $x \in \bR$ and observe that the marginal distribution of $\mu$ on $(\langle v, h_1 \rangle_H, \ldots, \langle v, h_N \rangle_H)$ is $(\bR^N, \prod_{k=1}^N \Psi_k(x_k)dx_k)$. By \eqref{A4'}, 
\begin{equation}\label{finite dimension divergence free}
\int_{\bR^N} e^{-2V^\dagger(\mathbf x)}\left[B^\dagger(\mathbf x) \cdot \nabla f(\mathbf x)\right]\prod_{k=1}^N \Psi_k(x_k)dx_k = 0, \ \forall f \in C_0^1(\bR^N). 
\end{equation}
Hence $e^{-2V^\dagger(\mathbf x)}\prod_{k=1}^N \Psi_k(x_k)dx_k$ is an (infinite) invariant measure for $X_t$. Noting that now $\{X_t^M; M > N\}$ forms a mutually independent system of Ornstein--Uhlenbeck processes which is independent of $X_t$, we can conclude that $\nu(dv) = e^{-2V(v)}\mu(dv)$ is a stationary measure for $\{u^v(t); t \geq 0\}$. 

For general $V$ and $B$, we consider their marginal expectations. Precisely speaking, for $N \geq 1$ and $\mathbf x = (x_1, \ldots, x_N) \in \bR^N$ define 
$$V_N^\dagger(\mathbf x) = \prod_{l=N+1}^\infty \left[\int_\bR \Psi_l(y_l)dy_l\right]V\left(\sum_{l=1}^N x_lh_l + \sum_{l=N+1}^\infty y_lh_l\right). $$
Similarly, write $U_k = e^{-2V}\langle B, h_k \rangle_H$ and for $k = 1, \ldots, N$ define 
$$B_{N, k}^\dagger(\mathbf x) = e^{2V_N^\dagger(\mathbf x)}\prod_{l=N+1}^\infty \left[\int_\bR \Psi_k(y_l)dy_l\right] U_k\left(\sum_{l=1}^N x_lh_l + \sum_{l=N+1}^\infty y_lh_l\right). $$
Define $V_N$ and $B_N$ by substituting $V_N^\dagger$, $B_N^\dagger$ for $V^\dagger$, $B^\dagger$ in \eqref{cylindrical V and B}, and let $u_N^v(t)$ be the solution to \eqref{spde} with $DV + B$ replaced by $DV_N + B_N$. From the Lipschitz continuity of $V$, $DV$ and $B$, we get for every fixed $v \in H$ that 
$$\lim_{N \rightarrow \infty} |V_N(v) - V(v)| + \|DV_N(v) - DV(v)\|_H + \|B_N(v) - B(v)\|_H = 0. $$
Since $V$, $DV$ and $B$ are uniformly bounded, we can obtain that 
$$\sup_{N \geq 1, v \in H}|V_N(v)| + \|DV_N(v)\|_H + \|B_N(v)\|_H < \infty. $$
Let $\mathbb E$ be the expectation with respect to the noise. With the above estimates, for a function $f$ on $E$ such that $f(v) = f^\dagger(\langle v, h_1 \rangle_H, \ldots, \langle v, h_M \rangle_H)$ with some $M \geq 1$ and $f^\dagger \in C_0^1(\bR^M)$, we can show that 
$$\lim_{N \rightarrow \infty} \mathbb E [f(u_N^v(t))] = \mathbb E [f(u^v(t))]. $$
Noting that $(V_N^\dagger, B_N^\dagger)$ fulfills \eqref{finite dimension divergence free}, the previous step implies that 
$$\int_E \mathbb E [f(u_N^v(t))]e^{-2V_N(v)}\mu(dv) = \int_E f(v)e^{-2V_N(v)}\mu(dv). $$
Because $f^\dagger$ is compactly supported, by taking $N \rightarrow \infty$ we get $\int_E \mathbb E [f(u^v(t))]\nu(dv) = \int_E f(v)\nu(dv)$. Applying Cauchy--Schwartz inequality, 
$$\int_E \left\{\mathbb E [f(u_N^v(t))]\right\}^2\nu(dv) \leq \int_E \mathbb E [f^2(u_N^v(t))]\nu(dv) = \int_E f^2(v)\nu(dv). $$
With this inequality and the fact that the finite-dimensional dependent functions is dense in $L^2(\nu)$, we can conclude that $\nu$ is invariant for $\{u^v(t); t \geq 0\}$. 

In case that $B = 0$, we can prove the reversibility of $\nu$ in the same way, only to observe that now the finite-dimensional approximation $X_t$ is a symmetric process, so that $e^{-2V^\dagger(\mathbf x)}\prod_{k=1}^N \Psi_k(x_k)dx_k$ becomes reversible. 
\end{proof}

The stationary measure obtained in Lemma \eqref{stationarity-frozen} is of infinite mass. To overcome this difficulty, we define a probability measure $\pi$ on $\Xi$ as 
\begin{equation}\label{invariant measure}
\pi(d\xi) = Z^{-1}\exp(-2\fV(\xi))\mu_0(dv_\xi) \otimes \bQ(d\xi(0)), 
\end{equation}
where $v_\xi \in E_0$ is determined by the one-to-one map from $\Xi$ to $\Sigma \times E_0$, and $Z$ is the normalization constant. Our aim is to show that $\pi$ is stationary for $\{\cP_t; t \geq 0\}$. 

\begin{prop}\label{stationarity}
$\pi$ is a stationary and ergodic probability measure for $\{\cP_t; t \geq 0\}$. Moreover, if $\fB \equiv 0$, then $\pi$ is reversible. 
\end{prop}

To prove Proposition \ref{stationarity}, we need to introduce the class of smooth functions on $\Xi$. Let $\cC$ be the dense subspace of $H$ defined as 
\begin{equation}\label{kernel functions}
\cC = \{h \in H~|~h \in C^2[0, 1], h'(0) = 0, h'(1) = 0, \langle h, \mathbf 1 \rangle_H \not= 0\}, 
\end{equation}
Pick arbitrarily some $\psi \in L^\infty(\Sigma; \bQ)$, $n \geq 1$, $h_1, \ldots, h_n \in \cC$ and $\ell^\dagger \in C_b^\infty(\bR^n)$ such that $\ell^\dagger$ together with all its derivatives belongs to $L^1(\bR^n; dx)$. Consider the function 
\begin{equation}\label{smooth functions}
\ff(\xi) = \int_\bR \psi(\tau_\theta\sigma) \ell(v - \theta\mathbf1)d\theta, \ \forall \xi = \tau_{[v(\cdot)]}\sigma \in \Xi, 
\end{equation}
where $\ell(v) = \ell^\dagger(\langle v, h_1 \rangle_H, \ldots, \langle v, h_n \rangle_H)$ for $v \in E$. Let $\cE_0(\Xi)$ be the collection of all $\ff$ in \eqref{smooth functions}, and call its linear span $\EXi$ \textit{the smooth function class} on $\Xi$. 

$\EXi$ is dense in $L^2(\Xi; \pi)$. To see that, pick $\fg \in L^2(\Xi; \pi)$ and suppose that $\langle \fg, \ff \rangle_\pi = 0$ for all $\ff \in \cE_0(\Xi)$. By virtue of \eqref{invariant measure} and \eqref{smooth functions}, it implies 
$$\frac{1}{Z}E_\bQ \left\{\int_{E_0} \mu_0(dv)\int_\bR e^{-2\fV(\tau_{[v(\cdot)]}\sigma)}\fg(\tau_{[v(\cdot)]}\sigma)\psi(\tau_\theta\sigma)\ell(v - \theta\mathbf1)d\theta\right\} = 0 $$
holds for all $\psi$ and $\ell$ satisfying the conditions above \eqref{smooth functions}. Applying the change of variable $\sigma = \tau_{-\theta}\sigma'$, we get from the stationarity of $\bQ$ under $\tau_\theta$ that 
$$\frac{1}{Z}E_\bQ \left\{\int_{E_0} \mu_0(dv)\int_\bR e^{-2\fV(\tau_{[v(\cdot) - \theta\mathbf1]}\sigma')}\fg(\tau_{[v(\cdot) - \theta\mathbf1]}\sigma')\psi(\sigma')\ell(v - \theta\mathbf1)d\theta\right\} = 0. $$
Since $v - \theta\mathbf1$ in this integral subjects to the infinite measure $\mu$ defined in Lemma \ref{stationarity-frozen}, we can rewrite the above equation as 
\begin{equation}\label{integral translation}
\frac{1}{Z}E_\bQ \left\{\int_E e^{-2\fV(\tau_{[v(\cdot)]}\sigma')}\fg(\tau_{[v(\cdot)]}\sigma')\psi(\sigma')\ell(v)\mu(dv)\right\} = 0. 
\end{equation}
From \eqref{integral translation} we obtain that $e^{-2\fV}\fg = 0$ for $\bQ$-almost all $\sigma$ and $\mu$-almost all $v$, and thus $\fg = 0$ in $L^2(\Xi, \pi)$. As $\EXi$ is the linear span of $\cE^0(\Xi)$, it is dense in $L^2(\Xi; \pi)$. 

\begin{proof}[Proof of Proposition \ref{stationarity}]
We begin with pointing out that to prove the stationarity of $\pi$, it suffices to show $E_\pi [\cP_t(\ff \cdot \ff')] = E_\pi [\ff \cdot \ff']$ for all $\ff$, $\ff' \in \EXi$. Indeed, by taking $\ff' = \ff$ we have $E_\pi [(\cP_t\ff)^2] \leq E_\pi [\cP_t(\ff^2)] = E_\pi [\ff^2]$, so $\cP_t$ is contractive on $L^2(\Xi; \pi)$. Since $\EXi$ is dense in $L^2(\Xi; \pi)$, it is easy to obtain $E_\pi [\cP_t\ff] = E_\pi \ff$ for all $\ff \in L^2(\Xi; \pi)$, thus $\pi$ is a stationary measure. 

Now we prove that $E_\pi [\cP_t(\ff \cdot \ff')] = E_\pi [\ff \cdot \ff']$ for $\ff$, $\ff' \in \EXi$. By the definition, it suffices to prove for $\ff$, $\ff' \in \cE_0(\Xi)$. Suppose that $\ff$, $\ff'$ are defined by \eqref{smooth functions} with some $(\psi, \ell)$ and $(\psi', \ell')$ respectively. Observe that for $\xi = \tau_{[v(\cdot)]}\sigma$, 
$$\ff \cdot \ff'(\xi) = \int_{\bR^2} \psi(\tau_\theta\sigma)\psi'(\tau_{\theta + \theta'}\sigma)\ell(v - \theta\mathbf1)\ell'(v - (\theta + \theta')\mathbf1)d\theta d\theta'. $$
To simplify the notations, in this proof we write $\Psi_\theta(\sigma) = \psi(\sigma)\psi'(\tau_\theta\sigma)$ and $L_\theta(v) = \ell(v)\ell'(v - \theta\mathbf1)$. With these notations we have 
$$E_\pi [\cP_t(\ff \cdot \ff')] = \frac{1}{Z}E_\bQ \left\{\int_{E_0} \mu_0(dv)\int_{\bR^2} e^{-2V(\sigma, v)}\Psi_{\theta'}(\tau_\theta\sigma)\mathbb E \left[L_{\theta'}(u^{\sigma, v}(t) - \theta\mathbf1)\right]d\theta d\theta'\right\}. $$
Since $u^{\sigma, v}(t) - \theta\mathbf1 = u^{\tau_\theta\sigma, v - \theta\mathbf1}(t)$, the change of variable $\sigma = \tau_{-\theta}\sigma'$ yields that
$$E_\pi [\cP_t(\ff \cdot \ff')] = \frac{1}{Z}E_\bQ \left\{\int_{E_0} \mu_0(dv) \int_{\bR^2} e^{-2V(\sigma', v - \theta\mathbf1)}\Psi_{\theta'}(\sigma')\mathbb E \left[L_{\theta'}(u^{\sigma', v - \theta\mathbf1}(t))\right]d\theta d\theta'\right\}. $$
With the same arguments as in calculating \eqref{integral translation} we get 
$$E_\pi [\cP_t(\ff \cdot \ff')] = \frac{1}{Z}E_\bQ \left\{\int_\bR d\theta'\left[\Psi_{\theta'}(\sigma)\int_E \mathbb E \left[L_{\theta'}(u^{\sigma, v}(t))\right]\pi^\sigma(dv)\right]\right\}. $$
Due to Lemma \ref{stationarity-frozen}, the last integral equals to $E_{\pi^\sigma} [L_{\theta'}]$, so that 
$$E_\pi [\cP_t(\ff \cdot \ff')] = \frac{1}{Z}E_\bQ \left\{\int_E \mu(dv)\int_\bR e^{-2V(\sigma, v)}\Psi_{\theta'}(\sigma)L_{\theta'}(v)d\theta'\right\}. $$
The desired equation $E_\pi [\cP_t(\ff \cdot \ff')] = E_\pi [\ff \cdot \ff']$ then follows. 

In case that $\fB = 0$, the reversibility of $\pi$ can be proved in the same way, only to notice that $\pi^\sigma$ is now reversible to the semigroup in frozen environment. 

To see the ergodicity, pick $G \subseteq \Xi$ such that $\pi(G) > 0$ and $\cP_t\mathbf1_G = \mathbf1_G$ for some $t > 0$, and we show that $\pi(G) = 1$. For every $\sigma \in \Sigma$, define $G_\sigma = \{v \in E; \tau_{[v(\cdot)]}\sigma \in G\}$. Take a strictly positive $f \in L^1(\bR; dx)$ and observe that 
$$\frac{1}{Z}\int_\bR \left\{E_\bQ \left[\int_{E_0} e^{-2V(\sigma, v)}\mathbf1_{G_\sigma}(v)\mu_0(dv)\right] \cdot f(c)\right\}dc = \pi(G) \cdot \int_\bR f(c)dc > 0. $$
Noting that $\mathbf1_{G_\sigma}(v) = \mathbf1_{G_{[\tau_c\sigma]}}(v - c\mathbf1)$, the above equation can be written as 
$$\int_\bR \left\{E_\bQ \left[\int_{E_0} e^{-2V(\tau_c\sigma, v - c\mathbf1)}\mathbf1_{G_{[\tau_c\sigma]}}(v - c\mathbf1)\mu_0(dv)\right] \cdot f(c)\right\}dc > 0. $$
By the $\tau_c$-invariance of $\bQ$ and the same strategy used in \eqref{integral translation} we get 
$$E_\bQ \left[\int_E e^{-2V(\sigma, v)}\mathbf1_{G_\sigma}(v)f(-v(0))\mu(dv)\right] > 0, $$
thus $\bQ(\sigma; \mu(G_\sigma) = 0) < 1$. The $\tau_c$-invariance of $\{\sigma; \mu(G_\sigma) = 0\}$ then implies that it is $\bQ$-null, thanks to the ergodicity. Since $\mu$ is absolutely continuous with respect to $\pi^\sigma$, we know further that $\bQ(\sigma; \pi^\sigma(G_\sigma) = 0) = 0$. Noting that for fixed $\sigma$, \eqref{spde'} satisfies the conditions mentioned in \cite[Sect. 1]{PZ95}, hence from \cite[Corollary 1.1 and Theorem 1.3]{PZ95}, $u^\sigma(t)$ is strong Feller and irreducible. Therefore, its transition probability $P^\sigma(t; v, \cdot)$ is equivalent to its stationary measure $\pi^\sigma$ for all $v$ and $t$, so that $\bQ(\sigma; P^\sigma(t;v, G_\sigma) = 0) = 0$. Meanwhile, by \eqref{semigroup}, 
$$0 = \langle \mathbf1_{G^c}, \cP_t\mathbf1_G \rangle_\pi = \frac{1}{Z}E_\bQ \left\{\int_{E_0} [1 - \mathbf1_{G_\sigma}(v)]e^{-2V(\sigma, v)}P^\sigma(t; v, G_\sigma)\mu_0(dv)\right\}. $$
As $e^{-2V(\sigma, v)}P^\sigma(t; v, G_\sigma)$ is strictly positive for every $v$, $t$ and almost all $\sigma$, the above equation implies that $\bQ(\sigma; \mu_0(G_\sigma) = 1) = 1$, and $\pi(G) = 1$ follows directly. 
\end{proof}

The last part of this section is devoted to the generator of $\cP_t$. Although $\cP_t$ is not strongly continuous under the topology of $C_b(\Xi)$ (cf. \cite{PT01}), due to the stationarity of $\pi$ we can extend $\cP_t$ to a $C_0$ semigroup of contractions on $L^2(\Xi; \pi)$. Denote the extension still by $\cP_t$ and define its generator $(\dom(\cK), \cK)$ through Hille--Yosida theorem. Recall that $\cK$ is unbounded and closed, and we say $A \subseteq \dom(\cK)$ is a core if $A$ is dense in $\dom(\cK)$ and $\cK$ coincides with the closure of $\cK|_A$. In the next proposition, we compute $\cK$ on the smooth function class on $\Xi$ and show that it forms a core. 

\begin{prop}\label{core}
$\EXi$ forms a core of $\cK$, and for all $\ff \in \EXi$, 
\begin{equation}\label{generator}
\cK\ff(\xi) = \frac{1}{2}\langle \partial_x^2[\cD\ff(\xi)], v_\xi \rangle_H - \langle \cD\ff(\xi), \cD\fV(\xi) + \fB(\xi) \rangle_H + \frac{1}{2} \tr[\cD^2\ff(\xi)]. 
\end{equation}
\end{prop}

Proposition \ref{core} is proved along the strategy in \cite{PT01}, where the maximal dissipativity for a class of Kolmogorov operators is discussed. We first prove the parallel result for a linear equation in Lemma \ref{core-linear}, and then extend it to the nonlinear case.

To discuss the linear case, consider an Ornstein--Uhlenbeck process $\{u_0^v(t, \cdot) \in E; t \geq 0\}$ with a given initial condition $v \in E$ satisfying that 
$$u_0^v(t) \triangleq S_tv + \int_0^t S_{t-r}dW_r, $$
where $W_r$ is the cylindrical Brownian motion appeared in \eqref{spde'}, and $\{S_t; t \geq 0\}$ is the semigroup on $H$ generated by $\frac{1}{2}\partial_x^2$ with Neumann boundary condition. Since $u_0^v(t)$ solves \eqref{spde'} when $\cD\fV + \fB \equiv 0$, it is a continuous $E$-valued Markov process. Denote by $\{\cT_t^0; t \geq 0\}$ the Markov semigroup determined by $u_0^v(t)$ on $\cB_b(E)$. In imitation of \eqref{semigroup}, for $\ff \in \cB_b(\Xi)$ and $t \geq 0$ we define 
$$\cT_t\ff(\xi) = \cT_t^0f^{\xi(0)}(v_\xi), \ \forall \xi \in \Xi. $$
In view of Proposition \ref{stationarity}, $\cT_t$ possesses a reversible measure $\pi_0(d\xi) = Z^{-1}\exp(-2\fV(\xi))\pi(d\xi)$. Write $\cH_{\pi_0} = L^2(\Xi; \pi_0)$, and extend $\{\cT_t; t \geq 0\}$ to a strongly continuous semigroup of contractions on $\cH_{\pi_0}$. Denote by $(\dom(\cL), \cL)$ its generator on $\cH_{\pi_0}$. Paralleling to Proposition \ref{core}, we have the next lemma for $\cL$ and $\pi_0$. 

\begin{lem}\label{core-linear}
$\EXi$ forms a core of $\cL$, and for all $\ff \in \EXi$, $\cL\ff(\xi) = \cL^0f^{\xi(0)}(v_\xi)$ where $\cL^0$ is the Ornstein--Uhlenbeck operator defined as 
$$\cL^0f(v) = \frac{1}{2}\left\langle \partial_x^2Df(v), v \right\rangle_H + \frac{1}{2}\tr\left[D^2f(v)\right]. $$
\end{lem}

\begin{proof}
First we present a basic estimate. Suppose that $f(v) = f^\dagger(\langle v, h_1 \rangle_H, \ldots, \langle v, h_n \rangle_H)$ for all $v \in E$, where $n \geq 1$, $f^\dagger \in C_b^\infty(\bR^n)$ and $h_1, \ldots, h_n \in \cC$ are chosen arbitrarily. Let $M = \max\{|\partial_\alpha f^\dagger|_\infty\}$, where the maximum is taking over all positive $n$-multiple indexes $\alpha$ such that $|\alpha| \leq 3$. For such $f$, we have (cf. \cite[Corollary 2.3]{PT01}) 
\begin{equation}\label{estimate-C0-linear}
\sup_{v \in E} \frac{1}{1 + \|v\|_H^2}\left|\frac{\cT_t^0f(v) - f(v)}{t} - \cL^0f(v)\right| \leq MC_t, 
\end{equation}
where $C_t$ is a constant depending on the $h_i$'s, satisfying that $C_t \rightarrow 0$ when $t \downarrow 0$ for fixed $h_1, \ldots, h_n \in \cC$. We here give the proof of \eqref{estimate-C0-linear} only for the case $n = 1$, since for larger $n$ the calculations are the same. When $n = 1$, we write $f(v) = f^\dagger(\langle v, h \rangle_H)$, and define $Q_{t, h} = \int_0^t \langle S_{t-r}h, dW(r, \cdot) \rangle_H$. With $\mathbb E$ being the expectation with respect to the noise, we have $\mathbb E [Q_{t, h}] = 0$ and 
\begin{equation*}
\begin{aligned}
&\left|\cT_t^0F(v) - F(v) - t\cL^0F(v)\right| \\
\leq\,&\left|\mathbb E \left[f(\langle v, S_th \rangle_H + Q_{t, h}) - f(\langle v, S_th \rangle_H)\right] - \frac{t}{2} \cdot f''(\langle v, h \rangle_H)\|h\|_H^2\right| \\
&+ \left|f(\langle v, S_th \rangle_H) - f(\langle v, h \rangle_H) - t \cdot f'(\langle v, h \rangle_H)\left\langle v, \frac{1}{2}h'' \right\rangle_H\right| \\
\leq\,&\frac{1}{2}\left|f''\right|_\infty \left|\mathbb E \left[Q_{t, h}^2\right] - t\|h\|_H^2\right| + \frac{1}{6}\left|f'''\right|_\infty \mathbb E\left[|Q_{t, h}|^3\right] \\
& + \left|f'\right|_\infty \|v\|_H \cdot \left\|S_th - h - \frac{t}{2}h''\right\|_H + \frac{1}{2}\left|f''\right|_\infty \|v\|_H^2 \cdot \left\|S_th - h\right\|_H^2. 
\end{aligned}
\end{equation*}
Since $h \in \dom(\frac{1}{2}\partial_x^2)$, \eqref{estimate-C0-linear} follows from direct calculation. 

Now pick some $\ff$ in the form of \eqref{smooth functions} and notice that for all $\sigma \in \Sigma$, 
$$f^\sigma(v) = \int_\bR \psi(\tau_\theta\sigma)\ell(v - \theta\mathbf1)d\theta, \ \forall v \in E. $$
In view of the conditions on $\psi$ and $\ell$, it is easy to observe that \eqref{estimate-C0-linear} holds for each $f^\sigma$ with common constants $M$ and $C_t$, hence 
$$\frac{1}{1 + \|v_\xi\|_H^2}\left|\frac{\cT_t\ff(\xi) - \ff(\xi)}{t} - \cL^0 f^{\xi(0)}(v_\xi)\right| \leq MC_t. $$
Applying dominated convergence theorem, we obtain that 
\begin{equation}\label{C0-linear}
\cH_{\pi_0}\text{-}\lim_{t \downarrow 0} \frac{\cT_t\ff - \ff}{t} = \cL^0 f^{\xi(0)}(v_\xi), \ \forall \ff \in \cE^0(\Xi). 
\end{equation}
As $\EXi$ is the linear span of $\cE^0(\Xi)$, \eqref{C0-linear} holds for all $\ff \in \EXi$, therefore $\EXi \subseteq \dom(\cL)$ and $\cL\ff(\xi) = \cL^0f^{\xi(0)}(v_\xi)$ for $\ff \in \EXi$. 

We are left to show that $\EXi$ is a core. By \cite[p. 17, Proposition 1.3.3]{EK05}, it suffices to show that $\EXi$ is preserved under $\{\cT_t; t \geq 0\}$. Notice that for $h \in \cC$, 
$$\langle u_0^v(t), h \rangle_H = \langle v, S_th \rangle_H + Q_{t,h} $$
and $S_th$ still belongs to $\cC$. Therefore, for $\ff$ in the form of \eqref{smooth functions}, 
$$\cT_t^0f^\sigma(v) = \mathbb E [f^\sigma(u_0^v(t))] = \int_\bR \psi(\tau_\theta\sigma)\mathbb E [\ell^\dagger(\vv_{t,h} - \theta\langle \vh \rangle_H + \vec Q_{t,h})]d\theta, $$
where $\vv_{t,h} = (\langle v, S_th_1 \rangle_H, \ldots, \langle v, S_th_n \rangle_H)$, $\langle \vh \rangle_H = (\langle h_1 \rangle_H, \ldots, \langle h_n \rangle_H)$ are deterministic $\bR^n$-vectors, and $\vec Q_{t,h} = (Q_{t,h_1}, \ldots, Q_{t,h_n})$ is a random vector. $\cT_t\ff \in \EXi$ then follows from the conditions on $\psi$ and $\ell^\dagger$, and consequently $\cT_t[\EXi] \subseteq \EXi$. 
\end{proof}

From Lemma \ref{core-linear}, $\cL(\ff^2) = 2\ff \cdot \cL\ff + \|\cD\ff\|_H^2$ holds for $\ff \in \EXi$, and the stationarity of $\pi_0$ implies that $E_{\pi_0} [\cL(\ff^2)] = 0$, so that Dirichlet form of $\cL$ reads 
\begin{equation}\label{dirichlet form-linear}
\langle \ff, -\cL\ff \rangle_{\pi_0} = \frac{1}{2}\int_\Xi \|\cD\ff(\xi)\|_H^2{\pi_0}(d\xi), \ \forall \ff \in \EXi. 
\end{equation}
Let $\Hpi = L^2(\Xi, \pi)$, and $\langle \cdot, \cdot \rangle_\pi$, $\|\cdot\|_\pi$ be the related inner product and norm. By \textbf{(A2)}, $\|\cdot\|_\pi$ and $\|\cdot\|_{\pi_0}$ are equivalent, so that elements in $\Hpi$ can be identified with those in $\cH_{\pi_0}$. With this observation, we close this section with the proof to Proposition \ref{core}. 

\begin{proof}[Proof of Proposition \ref{core}]
First recall the definition of differentiable function on $\Xi$ in Definition \ref{derivative}. Let $C_b^1(\Xi)$ denote the class of all differentiable functions $\ff$ on $\Xi$ such that $\cD\ff$ is continuous and uniformly bounded in $H$. Define 
$$\cK^\circ\ff = \cL\ff - \langle \cD\ff, \cD\fV + \fB \rangle_H, \ \forall \ff \in \dom(\cL) \cap C_b^1(\Xi). $$
Noticing that $\EXi \subseteq \dom(\cL) \cap C_b^1(\Xi)$ and that $\cK^\circ$ satisfies \eqref{generator} for $\ff \in \EXi$, to verify the explicit form \eqref{generator} for $\cK$ we only need to show that 
\begin{equation}\label{C0-whole}
\Hpi\text{-}\lim \frac{\cP_t\ff - \ff}{t} = \cK^\circ\ff, \ \forall \ff \in \dom(\cL) \cap C_b^1(\Xi). 
\end{equation}
Indeed, observing that $u^{\sigma, v}(t) - u_0^v(t) = \int_0^t S_{t-r}[-\cD\fV(\xi_r) - \fB(\xi_r)]dr$, for any $\ff \in C_b^1(\Xi)$ and fixed $\sigma$, $v$ we have 
\begin{equation*}
\begin{aligned}
&\lim_{t \downarrow 0} \frac{\cP_t^\sigma f^\sigma(v) - \cT_t^\sigma f^\sigma(v)}{t} \\
=\,&\lim_{t \downarrow 0} \mathbb E \left\langle Df^\sigma(u_0^v(t)), -\frac{1}{t}\int_0^t S_{t-r}[\cD\fV(\xi_r) + \fB(\xi_r)]dr \right\rangle_H + o(1) \\
=\,&\left\langle Df^\sigma(v), \cD\fV(\tau_{[v(\cdot)]}\sigma) + \fB(\tau_{[v(\cdot)]}\sigma) \right\rangle_H. \\
\end{aligned}
\end{equation*}
By \textbf{(A2)} and dominated convergence theorem, for $\ff \in C_b^1(\Xi)$ we obtain 
\begin{equation}\label{C0-nonlinear}
\Hpi\text{-}\lim_{t \downarrow 0} \frac{\cP_t\ff - \cT_t\ff}{t} = -\langle \cD\ff, \cD\fV + \fB \rangle_H. 
\end{equation}
As $\pi_0$ and $\pi$ are mutually equivalent, \eqref{C0-whole} follows from \eqref{C0-linear} and \eqref{C0-nonlinear} directly. 

Due to the Lumer--Phillips theorem (see, e.g., \cite{Y80} and \cite[p. 17, Proposition 1.3.1]{EK05}), to show that $\EXi$ is a core, it suffices to prove that $(\lambda - \cK)[\EXi]$ is dense in $\Hpi$ for some $\lambda > 0$. We first prove it under an additional condition. By \textbf{(A1)}, for $\bQ$-almost all $\sigma$, $DV(\sigma, \cdot) + B(\sigma, \cdot)$ can be extended to some $U^\sigma \in C_b(H; H)$ such that $\|U^\sigma\|_H$ is uniformly bounded in $\sigma$ and $v$. Assume further that \vspace{0.4em}\\
\textbf{(A1')} $U^\sigma \in C_b^2(H; H)$ and $\sup_{\Sigma \times H} |D^2U^\sigma(v)|_{L(H; H)} < \infty$. \vspace{0.4em}\\
Now for fixed $\sigma$, \eqref{spde'} becomes a stochastic heat equation with twice differentiable nonlinear term, thus obviously $\cP_t^\sigma [C_b^1(H)] \subseteq C_b^1(H)$. Furthermore, due to \cite[Theorem 1.2]{PZ95}, there is some constant $C$ such that 
$$\mathrm{sup}_H \|D[P_t^\sigma f]\|_H \leq C\max\{t^{-1/2}, 1\} \cdot \mathrm{sup}_H |f|, \ \forall f \in C_b^1(H). $$
Fix $\fg \in \EXi$, $\lambda > 0$ and define $\ff = (\lambda - \cK)^{-1}\fg$. We show that $\ff \in \dom(\cL) \cap C_b^1(\Xi)$. Notice the following formula holds for $\ff$ that 
\begin{equation}\label{resolvent}
\ff(\xi) = \int_0^\infty e^{-\lambda t} \cdot \cP_t\fg(\xi)dt, \ \forall \xi \in \Xi. 
\end{equation}
By the above estimates and dominated convergence theorem, $\ff \in C_b^1(\Xi)$ and $\|\cD\ff(\xi)\|_H$ is bounded from above by $C_\lambda\sup_\Xi |\ff|$ with some constant $C_\lambda$. Observe that 
\begin{equation}\label{decomposition-maximal dissipativity}
\frac{\cT_s \ff(\xi) - \ff(\xi)}{s} = \frac{\cP_s \ff(\xi) - \ff(\xi)}{s} + \frac{\cT_s \ff(\xi) - \cP_s\ff(\xi)}{s}. 
\end{equation}
Since $\ff \in \dom(\cK)$, the $\cH_{\pi_0}$-limit of the first term in the right-hand side of \eqref{decomposition-maximal dissipativity} exists when $s \downarrow 0$. Since $\ff \in C_b^1(\Xi)$, \eqref{C0-nonlinear} yields that the second term converges in $\Hpi$. Noting that $\cH_{\pi_0}$ and $\Hpi$ have equivalent norm, we have $\ff \in \dom(\cL)$. 

The previous proof now implies that $\cK\ff = \cL\ff - \langle \cD\ff, \cD\fV + \fB \rangle_H$. Since $\EXi$ is a core of $\cL$, we can pick $\ff_k \in \EXi$ such that $\ff_k \rightarrow \ff$ and $\cL\ff_k \rightarrow \cL\ff$, both in $\cH_{\pi_0}$. By \eqref{dirichlet form-linear}, $\|\cD(\ff_k - \ff)\|_H$ vanishes in $\cH_{\pi_0}$ and $\Hpi$. Therefore, 
$$E_\pi |\cK(\ff_k - \ff)|^2 \leq 2E_\pi |(\cL(\ff_k - \ff)|^2 + 2E_\pi |\langle \cD(\ff_k - \ff), \cD\fV + \fB \rangle_H|^2 \rightarrow 0 $$
as $k \rightarrow \infty$. Consequently, $\|\fg - (\lambda - \cK)\ff_k\|_\pi \rightarrow 0$ and $(\lambda - \cK)[\EXi]$ is dense in $\Hpi$. 

If \textbf{(A1')} fails to hold, we need to pick an approximating sequence. Indeed, suppose that we can find a sequence $(\fV_n, \fB_n)$ satisfying \textbf{(A1')}, \textbf{(A2)}, \textbf{(A4)} and that 
\begin{equation}\label{approximate condition-core}
\lim_{n \rightarrow \infty} E_\pi [\|\cD\fV + \fB - \cD\fV_n - \fB_n\|_H^2] = 0. 
\end{equation}
Let $\cK_n$ be the generator of the environment process related to the solution to \eqref{spde'} with $(\fV, \fB)$ replaced by $(\fV_n, \fB_n)$. Fix some $\fg \in \EXi$ and let $\ff_n = (\lambda - \cK_n)^{-1}\fg$. By the previous step, $\ff_n \in \dom(\cL) \cap C_b^1(\Xi)$, so that we have 
\begin{equation}\label{approximate formula-core}
(\lambda - \cK)\ff_n = \fg + (\cK_n - \cK)\ff_n = \fg + \langle \cD\fV_n + \fB_n - \cD\fV - \fB, \cD\ff_n \rangle_H. 
\end{equation}
Notice that for $\ff \in \dom(L) \cap C_b^1(\Xi)$ we have $\cK[\ff^2] = 2\ff \cdot \cK\ff + \|\cD\ff\|_H^2$. With the same way of computing the Dirichlet form for $\cL$ in \eqref{dirichlet form-linear}, we get that 
\begin{equation}\label{dirichlet form}
\langle \ff, -\cK\ff \rangle_\pi = 2^{-1}E_\pi \left[\|\cD\ff\|_H^2\right]. 
\end{equation}
In view of \eqref{approximate formula-core} and \eqref{dirichlet form}, 
$$\lambda \|\ff_n\|_\pi^2 + 2^{-1}E_\pi \left[\|\cD\ff_n\|_H^2\right] = \langle \fg, \ff_n \rangle_\pi + E_\pi [\langle \cD\fV_n + \fB_n - \cD\fV - \fB, \cD\ff_n \rangle_H \cdot \ff_n]. $$
By virtue of \eqref{resolvent}, $\sup_\Xi |\ff_n| \leq \lambda^{-1}\sup_\Xi |\fg|$ holds for all $n$, thus 
$$\frac{\lambda}{2}E_\pi \left[\|\cD\ff_n\|_H^2\right] \leq \sup_\Xi |\fg|^2 + \sup_\Xi |\fg| \cdot E_\pi [\langle \cD\fV_n + \fB_n - \cD\fV - \fB, \cD\ff_n \rangle_H]. $$
Using Cauchy--Schwarz inequality we obtain that 
$$\frac{\lambda}{4}E_\pi \left[\|\cD\ff_n\|_H^2\right] \leq \sup_\Xi |\fg|^2 + \frac{4}{\lambda}\sup_\Xi |\fg|^2 \cdot E_\pi \left[\|\cD\fV_n + \fB_n - \cD\fV - \fB\|_H^2\right]. $$
Hence, $\{E_\pi [\|\cD\ff_n\|_H^2]; n \geq 1\}$ is a uniformly bounded sequence. By \eqref{approximate formula-core} and \eqref{approximate condition-core}, $(\lambda - \cK)\ff_n$ converges to $\fg$ in the topology of $\Hpi$, thus $(\lambda - \cK)[\EXi]$ is dense. 

We are left to show the existence of $(\fV_n, \fB_n)$. Using the strategy in the proof of Lemma \ref{stationarity-frozen}, for fixed $\sigma$ without loss of generality we can assume $(V(\sigma, \cdot), B(\sigma, \cdot))$ to be in the form of \eqref{cylindrical V and B} with some $V^\dagger(\sigma, \cdot) \in C_b^1(\bR^N)$ and $B^\dagger(\sigma, \cdot) \in C_b(\bR^N; \bR^N)$. Pick smooth functions $V_n^\dagger(\sigma, \cdot)$ and $B_n^\dagger(\sigma, \cdot)$ such that $(\nabla V_n^\dagger + B_n^\dagger)(\sigma, \cdot)$ converges to $(\nabla V^\dagger + B^\dagger)(\sigma, \cdot)$ in a point-wise sense, $\sup_{\Sigma \times \bR^N}|\nabla V_n^\dagger + B_n^\dagger|$ is uniformly bounded from above and \eqref{finite dimension divergence free} holds with $V_n^\dagger(\sigma, \cdot)$ and $B_n^\dagger(\sigma, \cdot)$. Then, we get the desired sequence $(V_n(\sigma, \cdot), B_n(\sigma, \cdot))$ from $V_n^\dagger(\sigma, \cdot)$ and $B_n^\dagger(\sigma, \cdot)$ with \eqref{cylindrical V and B}. 
\end{proof}

\section{Central limit theorem for $u(t)/\sqrt t$}

Before giving the proof of Theorem \ref{clt}, we introduce two Hilbert spaces $\Hp$ and $\Hm$ related to the operator $\cK$ (cf. \cite[Sect. 2.2]{KLO12}). 

Notice that \eqref{dirichlet form} defines a semi-norm $\|\cdot\|_1^2$ on $\EXi$. Let $\sI_1 = \EXi/\sim_1$, where $\sim_1$ is the equivalence relation defined by $\ff \sim_1 \fg \iff \|\ff - \fg\|_1 = 0$. $\|\cdot\|_1$ naturally extends to a norm on $\sI_1$, and let $\Hp$ be the completion of $(\sI_1, \|\cdot\|_1)$. $\Hp$ becomes a Hilbert space under the inner product given by polarization: 
$$\langle \ff, \fg \rangle_1 \triangleq \frac{1}{4}(\|\ff + \fg\|_1^2 + \|\ff - \fg\|_1^2). $$
Noting that by \eqref{dirichlet form}, $\cD$ can be extended to a linear operator defined on $\Hp$ and taking values in $L^2(\Xi, H; \pi)$. We denote the extension still by $\cD$. Meanwhile, since for each $\ff \in \EXi$, $f^\sigma$ is finite-dimensional dependent, by applying the classical It\^o formula for finite-dimensional semi-martingales we get 
\begin{equation}\label{ito}
\ff(\xi_t) - \ff(\xi_0) = \int_0^t \cK\ff(\xi_r)dr + \int_0^t \langle \cD\ff(\xi_r), dW_r \rangle_H, \ \forall \ff \in \EXi. 
\end{equation}

The dual space $\Hm$ of $\Hp$ is defined as follows. For $\ff \in \Hpi$, let 
$$\|\ff\|_{-1}^2 \triangleq \sup_{\fg \in \EXi}\left\{2\langle \ff, \fg \rangle_\pi - \|\fg\|_1^2\right\}. $$
Let $\sI_{-1}^0 = \{\ff \in \Hpi, \|\ff\|_{-1} < \infty\}$ and $\sI_{-1} = \sI_{-1}^0/\sim_{-1}$, where $\sim_{-1}$ is the equivalence relation defined by $\ff \sim_{-1} \fg \iff \|\ff - \fg\|_{-1} = 0$. Let $\Hm$ be the completion of $(\sI_{-1}, \|\cdot\|_{-1})$, which is also a Hilbert space under the inner product $\langle \cdot, \cdot \rangle_{-1}$ defined through polarization. The following variational formula holds for $\|\cdot\|_{-1}$: 
\begin{equation}\label{variational formula for h-1}
\|\ff\|_{-1} = \sup_{\fg \in \EXi} \left\{\frac{\langle \ff, \fg \rangle_\pi}{\|\fg\|_1}; \|\fg\|_1 \not= 0\right\}, \ \forall \ff \in \Hpi\cap\Hm. 
\end{equation}

Recall that $\Sigma_0 = \{\sigma \in \Sigma~|~\tau_c\sigma = \sigma, \forall c \in \bR\}$. Observe that when $\sigma \notin \Sigma_0$, in view of Lemma \ref{key lemma}, $v_{\xi_t} = u^{\sigma, v}(t, \cdot) - u^{\sigma, v}(t, 0)$, hence the weak form of \eqref{spde'} reads 
\begin{equation} \label{decomposition-clt}
\langle u^\sigma(t), \varphi \rangle_H = \langle u^\sigma(0), \varphi \rangle_H + \int_0^t \fU_\varphi(\xi_r)dr + \langle W_t, \varphi \rangle_H, 
\end{equation}
where $\varphi \in C^2[0, 1]$ is a test function such that $\varphi'(0) = \varphi'(1) = 0$, and 
\begin{equation}\label{drift}
\fU_\varphi(\xi) = \frac{1}{2}\langle v_\xi, \partial_x^2\varphi \rangle_H - \langle \cD\fV(\xi) + \fB(\xi), \varphi \rangle_H. 
\end{equation}
When $\sigma \in \Sigma_0$ we have $v_{\xi_t} \equiv 0$ and \eqref{decomposition-clt} fails to hold. However, from \eqref{strong non-degenerate}, $\Sigma_0$ is a $\bQ$-null set here, therefore \eqref{decomposition-clt} holds $\bQ$-almost surely and this is sufficient for the $L^1(\bQ)$ convergence in Theorem \ref{clt}. 

Let $(\dom(\cK^*), \cK^*)$ be the adjoint operator of $\cK$ on $\Hpi$. By the standard arguments in \cite[Sects. 2.6, 2.7, 9.5]{KLO12}, to get a central limit theorem for $\langle u^\sigma(t), \varphi \rangle_H$ it is necessary to show that $\EXi$ is also a core of $\cK^*$, $\cK$ satisfies a sector condition and $\fU_\varphi \in \Hm$. We prove these results in this section. 

\begin{prop}\label{common core}
$\EXi$ forms a core of $\cK^*$, and for all $\ff \in \EXi$, 
\begin{equation}\label{adjoint}
\cK^*\ff(\xi) = \frac{1}{2}\langle \partial_x^2[\cD\ff(\xi)], v_\xi \rangle_H - \langle \cD\ff(\xi), \cD\fV(\xi) - \fB(\xi) \rangle_H + \frac{1}{2} \tr[\cD^2\ff(\xi)]. 
\end{equation}
\end{prop}

\begin{proof}
Let $\hat u^{\sigma, v}(t, \cdot)$ be the solution to \eqref{spde'} with $\fB$ replaced by $-\fB$, and $\{\hat \xi_t; t \geq 0\}$ be the corresponding environment process. By Proposition \ref{stationarity}, $\pi$ is stationary for $\hat \xi_t$, so that its Markov semigroup $\{\hat\cP_t; t \geq 0\}$ is a $C_0$ semigroup of contractions on $\Hpi$. Let $\hat\cK$ be its generator. By Proposition \ref{core}, $\hat\cK$ satisfies \eqref{adjoint} and $\EXi$ is a core of $\hat\cK$. 

It suffices to show $\cK^* = \hat\cK$. From Proposition \ref{stationarity}, $\pi$ is reversible if $\fB = 0$, thus 
\begin{equation}\label{symmetric part}
\langle \cK\ff + \langle \cD\ff, \fB \rangle_H, \fg \rangle_\pi = \langle \ff, \cK^*\fg - \langle \cD\fg, \fB \rangle_H \rangle_\pi, \ \forall \ff, \fg \in \EXi. 
\end{equation}
Meanwhile, in view of \textbf{(A4)}, for $\ff$ and $\fg \in \EXi$ we have 
\begin{equation}\label{anti-symmetric part}
E_\pi \left[\langle \cD [\ff\fg], \fB \rangle_H\right] = \frac{1}{2Z}E_\bQ \left\{\int_{E_0} e^{-2V(\sigma, v)}\langle D[f^\sigma g^\sigma](v), B(\sigma, v) \rangle_H\mu_0(dv)\right\} = 0. 
\end{equation}
Combing \eqref{symmetric part} and \eqref{anti-symmetric part} shows that $\hat\cK = \cK^*$ on $\EXi$. Noting that $\EXi$ is a core of $\hat\cK$ and $\cK^*$ is always closed, this would imply that $\hat\cK \subseteq \cK^*$. On the other hand, since $\cK^*$ is the adjoint of a Markov generator of a $C_0$ semigroup contractions, it is easy to show that $\lambda - \cK^*$ is invertible for all $\lambda > 0$. For all $\ff \in \Hpi$, 
$$\ff = (\lambda - \hat\cK)\left[(\lambda - \hat\cK)^{-1}\ff\right] = (\lambda - \cK^*)\left[(\lambda - \hat\cK)^{-1}\ff\right], $$
therefore $(\lambda - \cK^*)^{-1} = (\lambda - \hat\cK)^{-1}$, and consequently $\dom(\cK^*) = \dom(\hat\cK)$. 
\end{proof}

\begin{prop}\label{sector condition}
The generator $\cK$ satisfies the sector condition, i.e., 
$$\langle \cK\ff, \fg \rangle_\pi^2 \leq C\langle -\cK\ff, \ff \rangle_\pi \langle -\cK\fg, \fg \rangle_\pi, \ \forall \ff, \fg \in \dom(\cK), $$
where $C = C(\fV, \fB)$ is a finite constant depending only on $\fV$ and $\fB$. 
\end{prop}

\begin{proof}
For $\ff \in \dom(\cK) \cap \dom(\cK^*)$, let $\cK_s\ff = \frac{1}{2}(\cK\ff + \cK^*\ff)$ and $\cK_a\ff = \frac{1}{2}(\cK\ff - \cK^*\ff)$. By Proposition \ref{common core} it suffices to show for $\ff$, $\fg \in \EXi$ that 
$$\langle \cK_a\ff, \fg \rangle_\pi^2 \leq C(\fV, \fB)E_\pi\|\cD\ff\|_H^2E_\pi\|\cD\fg\|_H^2. $$
For a function $g$ on $E$, let $\langle g \rangle_{\mu_0} \triangleq E_{\mu_0} [g]$. By virtue of \textbf{(A4)}, 
$$\langle -\cK_a\ff, \fg \rangle_\pi = \frac{1}{Z}E_\bQ \left[\int_{E_0}  e^{-2V(\sigma, v)}\langle Df^\sigma(v), B(\sigma, v) \rangle_H\left(g^\sigma(v) - \langle g^\sigma \rangle_{\mu_0}\right)\mu_0(dv)\right]. $$
By Cauchy--Schwarz inequality and \textbf{(A2)} we obtain that 
\begin{equation}\label{estimate for sector condition}
\langle -\cK_a\ff, \fg \rangle_\pi \leq C_1\|\ff\|_1 \left\{\frac{1}{Z}E_\bQ \left[\int_{E_0} \left|g^\sigma(v) - \langle g^\sigma \rangle_{\mu_0}\right|^2\mu_0(dv)\right]\right\}^{1/2} 
\end{equation}
with $C_1 = 2\sup \|\fB\|_H\exp(\sup |\fV|)$. 

Now we apply the Poincar\'e inequality for Wiener measure. Let $\hat H$ be the Cameron--Martin space $\hat H = \{h \in W^{1,2}([0, 1]; \bR)~|~h(0) = 0\}$, equipped with the Cameron--Martin norm $\|h\|_{\hat H} \triangleq \|\dot h\|_H$, where $\dot h$ stands for the weak derivative of $h$. For a function $g$ on $E_0$ having the form $g(v) = g^\dagger(\langle v, h_1 \rangle_H, \ldots, \langle v, h_M \rangle_H)$ for some $M \geq 1$, $g^\dagger \in C_0^1(\bR^M)$ and $h_i \in H$, recall its $\hat H$-derivative $D_{\hat H}g \in \hat H$ is 
$$D_{\hat H}g(v) = \sum_{i=1}^M \left[\partial_ig^\dagger(\langle v, h_1 \rangle_H, \ldots, \langle v, h_N \rangle_H)\int_0^\cdot \hat h_i(x)dx\right], \ \forall v \in E_0, $$
where $\hat h_i = \int_\cdot^1 h_i(x)dx \in \hat H$. The Poincar\'e inequality for Wiener measure (see, e.g., \cite[p. 226, Theorem 5.5.1]{B98}) yields that 
\begin{equation}\label{poincare inequality}
\int_{E_0} |g(v) - \langle g \rangle_{\mu_0}|^2\mu_0(dv) \leq \int_{E_0} \|D_{\hat H}g(v)\|_{\hat H}^2\mu_0(dv). 
\end{equation}
To continue, write $\vv_h = (\langle v, h_1 \rangle_H, \ldots, \langle v, h_M \rangle_H)$ for shot. Observe that from the explicit formula for $D_{\hat H}g(v)$ and Cauchy--Schwarz inequality, 
\begin{equation*}
\begin{aligned}
\|D_{\hat H}g(v)\|_{\hat H}^2 &= \int_0^1\left[\sum_{i=1}^M \partial_ig^\dagger(\vv_h)\int_y^1 h_i(x)dx\right]^2dy \\
&\leq \int_0^1 (1 - y)\int_y^1\left[\sum_{i=1}^M \partial_ig^\dagger(\vv_h)h_i(x)\right]^2dxdy \\
&\leq \int_0^1 \left[\sum_{i=1}^M \partial_ig^\dagger(\vv_h)h_i(x)\right]^2dx = \|Dg(v)\|_H^2. 
\end{aligned}
\end{equation*}
Hence in the right-hand side of \eqref{poincare inequality}, we can substitute $\|D_{\hat H}g(v)\|_{\hat H}^2$ with $\|Dg(v)\|_H^2$. As $\fg \in \EXi$, by applying \eqref{poincare inequality} to $g^\sigma$ in \eqref{estimate for sector condition} we conclude that 
$$\langle -\cK_a\ff, \fg \rangle_\pi \leq C_1\|\ff\|_1\|\fg\|_1, $$
and the sector condition then follows. 
\end{proof}

\begin{prop}
For every test function $\varphi \in C^2[0, 1]$ which satisfies that $\varphi'(0) = 0$, $\varphi'(1) = 0$, we have $\fU_\varphi$ defined by \eqref{drift} belongs to $\Hm$. 
\end{prop}

\begin{proof}
In view of \eqref{variational formula for h-1}, we only need to show that there exists some finite constant $C$ such that $\langle \fU_\varphi, \fg \rangle_\pi \leq C\|\fg\|_1$ holds for all $\fg \in \EXi$. Taking $f(v) = \langle v, \varphi \rangle_H$ in \textbf{(A4)}, we obtain that for $\bQ$-almost all $\sigma$, 
\begin{equation}\label{A4''}
E_{\mu_0} \left[e^{-2V(\sigma, v)}\langle B(\sigma, v), \varphi \rangle_H\right] = 0. 
\end{equation}
Pick $\fg \in \EXi$ and \eqref{A4''} implies that 
$$E_\pi [\langle \fB, \varphi \rangle_H \cdot \fg] = \frac{1}{Z}E_\bQ \left[\int_{E_0} e^{-2V(\sigma, v)}\langle B(\sigma, v), \varphi \rangle_H(g^\sigma(v) - \langle g^\sigma \rangle)\mu_0(dv)\right]. $$
With the same argument as in the proof of Proposition \ref{sector condition}, we obtain that 
\begin{equation}\label{h-1 estimate for anti-symmetric drift}
E_\pi [\langle \fB, \varphi \rangle_H \cdot \fg] \leq C\|\varphi\|_H\|\fg\|_1. 
\end{equation}
On the other hand, by the integral-by-part formula for Wiener measure (see, e.g., \cite[p. 208, Proposition 5.1.6]{B98}) it is easy to get 
\begin{equation}\label{h-1 estimate for symmetric drift}
E_\pi [(\fU_\varphi - \langle \fB, \varphi \rangle_H) \cdot \fg] = \frac{1}{2}E_\pi [\langle \varphi, \cD\fg \rangle_H] \leq \|\varphi\|_H\|\fg\|_1. 
\end{equation}
Combining \eqref{h-1 estimate for anti-symmetric drift} and \eqref{h-1 estimate for symmetric drift}, we get the desired estimate. 
\end{proof}

We are at the position to prove Theorem \ref{clt}. We first prove the weak convergence of $\langle u^\sigma(t), \varphi \rangle_H/\sqrt t$ for fixed $\varphi$. For $\lambda > 0$, consider the resolvent equation: 
\begin{equation}\label{resolvent equation-clt}
\lambda \ff_{\lambda, \varphi} - \cK\ff_{\lambda, \varphi} = \fU_\varphi. 
\end{equation}
Define the Dynkin's martingale corresponding to $\ff_{\lambda, \varphi} \in \dom(\cK)$ as 
$$M_{\lambda, \varphi}(t) = \ff_{\lambda, \varphi}(\xi_t) - \ff_{\lambda, \varphi}(\xi_0) - \int_0^t \cK\ff_{\lambda, \varphi}(\xi_r)dr. $$
$M_{\lambda, \varphi}(t)$ is a square integrable martingale. Due to the arguments in \cite[Theorem 2.14, Sects. 2.7, 9.5]{KLO12}, with Proposition \ref{sector condition} we can obtain that 
\begin{equation} \label{convergence of M}
\lim_{\lambda \downarrow 0} \mathbb{E}^\pi \left[\sup_{0 \leq t \leq T}\left|M_{\lambda, \varphi}(t) - \int_0^t \langle \cD\ff_\varphi(\xi_r), dW_r \rangle_H\right|^2\right] = 0, \ \forall T > 0, 
\end{equation}
where $\ff_\varphi = \Hp\text{-}\lim_{\lambda \downarrow 0}\ff_{\lambda, \varphi}$. Furthermore, by \cite[p. 51, Lemma 2.10]{KLO12}, 
\begin{equation} \label{vanishment of R}
\lim_{t \rightarrow \infty} \lim_{\lambda \downarrow 0}\mathbb{E}^\pi \left[\frac{1}{t}\left|\int_0^t \fU_\varphi(\xi_r)dr - M_{\lambda, \varphi}(t)\right|^2\right] = 0. 
\end{equation}
Combining \eqref{decomposition-clt}, \eqref{convergence of M} and \eqref{vanishment of R} together, we get 
$$\lim_{t \rightarrow \infty}\mathbb{E}^\pi \left[\frac{1}{t}\left|\langle u^{\sigma, v}(t), \varphi \rangle_H - \int_0^t \langle \cD\ff_\varphi(\xi_r) + \varphi, dW_r \rangle_H\right|^2\right] = 0. $$
Notice that $\pi$ is ergodic, by central limit theorem for martingales (see, e.g., \cite[p. 36, Theorem 2.1]{KLO12}), for all bounded continuous function $f$ on $\bR$, 
\begin{equation} \label{clt-direction}
\lim_{t \rightarrow \infty}E_{\bQ} \left|E_P \left[f\left(\frac{\langle u(t), \varphi \rangle_H}{\sqrt{t}}\right)\right] - \int_{\mathbb{R}} f(y)\Phi_{a(\varphi)}(y)dy\right| = 0, 
\end{equation}
where $\Phi_{a(\varphi)}$ stands for the probability density function of Gaussian distribution with mean 0 and covariance $a^2(\varphi) = E_\pi \|\cD\ff_\varphi + \varphi\|_H^2$. 

Next, to prove \eqref{clt-equation}, it suffices to verify that $a(\varphi) = 0$ for all $\varphi$ orthogonal to the constant function $\mathbf 1$. Consider the function $\fg_\varphi(\xi) \triangleq -\langle v_\xi, \varphi \rangle_H$ for such $\varphi$. Noting that $\fg_\varphi$ is well defined and solves the cell problem $-\cK\fg_\varphi = \fU_\varphi$, hence $\fg_\varphi = \ff_\varphi$ in \eqref{convergence of M}. Consequently, $a^2(\varphi) = E_\pi \|\cD\fg_\varphi + \varphi\|_H^2 = 0$ and thus \eqref{clt-equation} holds with 
\begin{equation} \label{diffusion constant}
a^2 = a^2(\mathbf 1) = \lim_{\lambda \downarrow 0} E_\pi \|\cD\ff_{\lambda, \mathbf 1} + 1\|_H^2, 
\end{equation}
where $\ff_{\lambda, \mathbf 1}$ is the solution to \eqref{resolvent equation-clt} when $\varphi = \mathbf 1$. 

We are left to show that $a^2 \in [C, 1]$ for some $C > 0$ depending only on $\fV$. To this end, we introduce two scalar products for $H$-valued functions on $\Xi$ (cf. \cite[Sect. 2.2]{KLO12}). Consider $F$, $G: \Xi \rightarrow H$. Let $\llangle F, G \rrangle_\pi = E_\pi [\langle F, G \rangle_H]$, and 
$$\llangle F, G \rrangle = E_\bQ \left[\int_\Xi \left\langle F(\tau_{[v(\cdot)]}\sigma), G(\tau_{[v(\cdot)]}\sigma) \right\rangle_H \mu_0(dv_\xi) \right]. $$
Furthermore, define the norm $|\!|\!|F|\!|\!|^2 = \llangle F, F \rrangle$ and denote by $L^2(\Xi, H)$ the Hilbert space consisting of all functions with finite $|\!|\!|\cdot|\!|\!|$ norm, equipped with the inner product $\llangle \cdot, \cdot \rrangle$. Define $|\!|\!|\cdot|\!|\!|_\pi$ and $L_\pi^2(\Xi, H)$ from $\llangle \cdot, \cdot \rrangle_\pi$ similarly. 

From \eqref{resolvent equation-clt} we know $\langle \lambda \ff_{\lambda, \varphi} - \cK \ff_{\lambda, \varphi} - \fU_\varphi, \ff_{\lambda, \varphi} \rangle_\pi = 0$, which implies that 
$$\lambda \|\ff_{\lambda, \varphi}\|_\pi^2 + \frac{1}{2}|\!|\!| \cD\ff_{\lambda, \varphi} |\!|\!|_\pi^2 + \frac{1}{2}\llangle \varphi, \cD\ff_{\lambda, \varphi} \rrangle_\pi = 0. $$
Let $\lambda \downarrow 0$, we have $\llangle \cD\ff_\varphi + \varphi, \cD\ff_\varphi \rrangle_\pi = 0$, so that $\cD\ff_\varphi + \varphi$ and $-\cD\ff_\varphi$ are orthogonal to each other in $L_\pi^2(\Xi, H)$. Taking $\varphi = \mathbf 1$, it easily follows from the orthogonality that 
$$a^2 = E_\pi \|\cD\ff_\varphi + \varphi\|_H^2 \leq E_\pi \|\cD\ff_\varphi + \varphi - \cD\ff_\varphi\|_H^2 = 1. $$
On the other hand, since $\bQ$ is $\{\tau_c\}$-invariant, for all $\ff \in \EXi$ we have 
$$\llangle \mathbf 1, \cD\ff \rrangle = \int_\Sigma\int_{E_0} \left[\left.\!\frac{d}{dc} f^{\tau_c\sigma}(v)\right|_{c = 0}\right]\mu_0(dv)\bQ(d\sigma) = 0. $$
Since $\EXi$ is dense in $\Hp$, $\mathbf 1$ is orthogonal to the subspace $\cH_g \triangleq \{\cD\ff; \ff \in \Hp\}$ under the inner product $\llangle \cdot, \cdot \rrangle$. Therefore, with $C = Z^{-1} \cdot \exp(-2\sup|\fV|)$ we have 
$$a^2 \geq \min_{\ff \in \Hp} |\!|\!|\cD\ff + \mathbf 1|\!|\!|_\pi^2 \geq C\min_{\ff \in \Hp} |\!|\!|\cD\ff + \mathbf 1|\!|\!|^2 = C. $$
The proof of Theorem \ref{clt} is then completed. 

\begin{rem}
In the symmetric case $\fB = 0$, we can prove the following variational formula (cf. \cite[p. 335, Theorem 10.5]{KLO12}) for the limiting variance $a^2(\varphi)$: 
\begin{equation}\label{variational principle}
a^2(\varphi) = \min_{\ff \in \Hp}|\!|\!|\cD\ff + \varphi|\!|\!|_\pi^2. 
\end{equation}
Indeed, let $\fB = 0$ and write \eqref{resolvent equation-clt} in the weak form for $\fg \in \EH$, 
$$\lambda\langle \ff_{\lambda, \varphi}, \fg \rangle_\pi + \frac{1}{2}\llangle \cD\ff_{\lambda, \varphi} + \varphi, \cD\fg \rrangle_\pi = 0. $$
Let $\lambda \downarrow 0$, we conclude that $\cD\ff_\varphi + \varphi$ is orthogonal to $\cH_g$ under the inner product $\llangle \cdot, \cdot \rrangle_\pi$. Hence, \eqref{variational principle} is clear because that $\cD\ff_\lambda \in \cH_g$. In summary, the non-gradient term in \eqref{spde'} always enhances the fluctuation of the solution. 
\end{rem}

\section*{Acknowledgements}
The author greatly thanks Professor Tadahisa Funaki and Professor Stefano Olla for their instructive discussion and suggestions. The author also thanks Professor Jean-Dominique Deuschel for his comments on quenched results. 

\bibliography{[BIB]clt_for_spde.bib}

\end{document}